\hsize=16true cm \vsize=23true cm
%
\font\peque=cmr8
\font\sc=cmcsc10 \font\msbmnormal=msbm10
\font\msbmpeq=msbm7 \font\msbmmuypeq=msbm5
\newfam\numeros
\textfont\numeros=\msbmnormal \scriptfont\numeros=\msbmpeq
\scriptscriptfont\numeros=\msbmmuypeq
\def\num{\fam\numeros\msbmnormal}

\def\N{{\num N}}
\def\R{{\num R}}

\def\FM{{\cal  M}}

\def\FA{{\cal  A}}
\def\FP{{\cal P}}
%
\font\tenBb=msbm10 \font\sevenBb=msbm7 \font\fiveBb=msbm5
\newfam\Bbfam \textfont\Bbfam=\tenBb \scriptfont\Bbfam=\sevenBb
\scriptscriptfont\Bbfam=\fiveBb
 2 \font\grandebf=cmbx8 scaled
\magstep 2
\def\titulo#1{\ \smallskip
               \centerline{\grandebf #1}
               \smallskip}
\def\autor#1{\ \bigskip
              \centerline{\sc #1}
               \bigskip}
\def\parrafo#1{\
              \smallskip
              {\sc #1}
              \medskip}
%
\def\dem{{\it Proof:\/ \ }}
%
\def\bull{\hfill\vrule height .9ex width .8ex depth -.1ex } 

\def\c00{c_{00}}
\def\norma#1{\left\Vert #1 \right\Vert}
\def\vabs#1{\left\vert #1 \right\vert}
\def\3norma#1{\left\vert\hskip-1pt {\norma {#1}}\hskip-1pt \right\vert}

\let\lra=\longrightarrow
\let\ssi=\Longleftrightarrow
\let\ds=\displaystyle
\def\sop{\mathop{\rm supp}}

\titulo {ON TOTAL INCOMPARABILITY OF}
\titulo {MIXED TSIRELSON SPACES}
\autor {{Julio Bernu\'es\footnote{*}{\peque Partially supported by
DGES grant (Spain)} and Javier Pascual.}} \centerline{ University
of Zaragoza (Spain)}

\parrafo {ABSTRACT}

We give criteria of total incomparability for certain classes of
 mixed Tsirelson spaces. We show that spaces of the
form $T[(\FM_k,\theta_k)_{k =1}^{\ell}]$ with index $i(\FM_k)$
finite are either $c_0$ or $\ell_p$ saturated for some $p$ and we
characterize when any two spaces of such a form are totally
incomparable in terms of the index $i(\FM_k)$ and the parameter
$\theta_k$. Also, we give sufficient conditions of total
incomparability for a particular class of spaces of the form
$T[(\FA_k,\theta_k)_{k = 1}^\infty]$ in terms of the asymptotic
behaviour of the sequence $\Vert\sum_{i=1}^n e_i\Vert$ where
$(e_i)$ is the canonical basis.
\bigskip

{\it Key words and phrases,} Mixed Tsirelson spaces, totally
incomparable spaces.\bigskip

 \bigskip

\parrafo {0. Introduction}

Denote by $\c00$ the vector space of all real valued sequences
which are eventually zero and by $(e_i)_{i=1}^\infty$ its usual
unit vector basis. For $E \subset \N$ and $\ds x=\sum_{i=1}^\infty
a_i e_i\in\c00$ we denote $\ds Ex = \sum_{i \in E} a_ie_i$. Also,
for finite subsets $E,F \subseteq \N$, we write $E<F$ (or $E\le
F$) if $\max E < \min E$ ($\max E \le \min E$). For simplicity, we
write $n\le E$ instead of $\{n\}\le E$.

 Mixed Tsirelson spaces were introduced in full
generality in [2]. We can define those spaces, denoted by
$T[(\FM_k,\theta_k)_{k \in I}]$, as the completion of $\c00$ under
a norm which satisfies an implicit equation of the following kind:
$$ \norma x =\max \left\{ \norma x_\infty, \sup_{k \in I}
\left\{\theta_k \sup_{n \in \N} \left\{ \sum_{i=1}^n
\norma{E_ix}\vert\ (E_i)_{i=1}^n\ \FM_k-admissible \right\}
\right\} \right\},\qquad x\in\c00$$ where the $\FM_k$'s are
certain (see Definition 4 below) families of finite subsets of
$\N$, $\theta_k \in (0,1]$ for all $k \in I\subseteq\N$ and
$(E_i)_{i=1}^n$ is $\FM_k$-admissible if there exists $\{ m_1,
\dots,m_n\} \in \FM_k$ such that $m_1 \leq E_1< m_2 \leq E_2 <
\dots < m_n \leq E_n$.

The first remarkable space in this class is the so called
Tsirelson space, introduced by Figiel and Johnson [7] in 1974.
(It is actually the dual of the space originally constructed by
Tsirelson in [12].) In our notation this space is $T[{\cal
S},1/2]$, where ${\cal S}$ is Schreier's class, that is, the set
of subsets of $\N$ of cardinality smaller than their first
element.  Since its construction it was usually considered a
\lq\lq pathological space", a place to look for counterexamples to
statements in the Banach space theory. In fact, the reason why it
was constructed was to provide a counterexample to the assertion
\lq\lq every Banach space contains $c_0$ or $\ell_p$ for some
$1\le p<\infty$".

The second space of the class is Tzafriri space, introduced in
1979 in [13] ($T[({\cal A}_k,\gamma/\sqrt k)_{k \in\N}],
0<\gamma<1$ in our notation where ${\cal A}_k$ is the set of
subsets of $\N$ of at most $k$ elements), also constructed as a
counterexample to a statement in the Banach space theory. In 1991
a third example, namely the Schlumprecht space $T[({\cal
A}_k,1/\log_2(1+k))_{k \in\N}]$, was considered, see [11], and
with its help  a fruitful period started when many \lq\lq
classical" problems in the infinite dimensional Banach space
theory were solved, such as the distortion problem or the
unconditional basic sequence problem.

A common feature of the three Banach spaces mentioned above is
that they do not contain any $\ell_p, 1\le p<\infty$ or $c_0$.
(Actually, in the case of Tzafriri spaces this has been proved,
as far as we know, only for $0<\gamma<10^{-6}$, see [6].)
Moreover, since $\ell_p, 1\le p<\infty$ and $c_0$ are minimal
(recall that a Banach space $X$ is minimal if every subspace of
$X$ contains a further subspace isomorphic to $X$) it easily
follows that they are totally incomparable to any of the three
examples above (recall that two Banach spaces $X$ and $Y$ are
totally incomparable if no subspace of $X$ is isomorphic to any
of $Y$). We use the word \lq\lq subspace" here and throughout the
paper for \lq\lq closed infinite dimensional subspace".

In 1986 Bellenot [3] showed that $\ell_p, 1\le p<\infty$ and
$c_0$ are isomorphic to mixed Tsirelson spaces of the form
$T[({\cal A}_n,\theta)], \theta\in(0,1]$. This was somewhat
surprising as it showed that $\ell_p, 1\le p<\infty$ and $c_0$
belong to a class of spaces up to then considered pathological.

It is well known that $\ell_p, 1\le p<\infty$ and $c_0$ are
totally incomparable to each other. Moreover, $\ell_p$ and $c_0$
and the three examples, with $0<\gamma<10^{-6}$ in the case of
Tzafriri space, are all
 totally incomparable to each
other (see [6] for the details and also use the minimality of the
Schlumprecht space). This shows that, at least in the examples
considered, the modification of the $\theta_k's$ or the $ \FM_k's$
produce totally incomparable spaces.

In the first section we discuss in full generality the case when
$\theta_k=1$ for some $k$. In this case, the spaces $c_0$ and
$\ell_1$ will play a crucial role.

In the second section we consider mixed Tsirelson spaces of the
form $T[(\FM_k,\theta_k)_{k =1}^{\ell}], \theta_k \in (0,1),$
with index $i(\FM_k)$, as defined in [2], finite and we
characterize when any two spaces of such a form are totally
incomparable. This is done by following the ideas in [4] and
showing that every such space is either $c_0$ or $\ell_p$
saturated for some $p$. Recall that given a Banach space $Y$, a
Banach space $X$ is $Y$ saturated if every subspace of $X$
contains a further subspace isomorphic to $Y$.

In the third section we focus on spaces of the form
$T[(\FA_k,\theta_k)_{k = 1}^\infty], \theta_k \in (0,1],$ such
that $\ell_1$ is finitely block represented in every block
subspace. We give sufficient conditions of total incomparability
in terms of the asymptotic behaviour of the sequence
$\Vert\sum_{i=1}^n e_i\Vert$ where $(e_i)$ is the canonical
basis. These conditions apply to
 cases different from those considered in [9].\bigskip

\noindent{\bf Notation}. If $K$ is a subset of a Banach space $X$,
$\overline {Span} \{ K \}$ denotes the closure of the algebraic
linear span of $K$.  If $x = \sum_{i=1}^\infty a_ie_i \in \c00$,
the support of $x$ is the set $\sop(x) = \{ i \in \N\ \vert\ a_i
\not= 0 \}$.  For $x,y \in \c00$ we write $x < y$ if $\sop (x) <
\sop (y)$. We say that $E_1,\dots,E_n\subset\N$ are successive if
$E_1<E_2<\dots <E_n$. The vectors $x_1,\dots,x_n$ are successive
if their supports are. A block sequence $(x_i)$ is a sequence of
successive vectors.  The cardinality of a set $E$ is denoted by
$\vabs E$. The standard norm of $\ell_p, 1\leq p \le\infty$ is
denoted by $\norma\cdot_p$. Other unexplained notation is standard
and can be found for instance in [8].

\proclaim Definition 1. Let $\FM$ be a family of finite subsets
of $\N$. We say that $\FM$ is compact if the set $\{ \aleph_E \
\vert \ E \in \FM \}$ is a compact subset of the Cantor set $\{
0,1\}^\N$ with the product topology.

\proclaim Remark 1. In Definition 1, $\{ 0,1\}^\N$ is identified
with the space of all mappings $f:\N \lra \{0,1\}$ and $\aleph_E$
is the characteristic function of $E$. In $\{ 0,1\}^\N$, the
convergence under the product topology is the pointwise
 convergence. Therefore if $E \subseteq \N$ is a finite set and
$\aleph_{E_k} $ converges to $\aleph_E$ pointwise, there exists $N
\in \N$ such that $E \subseteq E_k$ for all $k \ge N$.

\proclaim Definition 2. Let $\FM $ be a family of finite subsets
of $\N$. We say that $\FM$ is hereditary if $E \in \FM $ and $F
\subseteq E$ implies that $F \in \FM$.

\proclaim Definition 3. Let $\FM$ be a compact family of finite
subsets of\ $\N$. We define a transfinite sequence\break
$(\FM^{(\lambda)})$  of subsets of $\FM$ as follows:
\item{1.} $\FM^{(0)}=\FM$.
\item{2.} $\FM^{(\lambda+1)}=\{ E \in \FM \ \vert\ \aleph_E$ is a limit point
                    of the set $\{\aleph_E \vert E \in
                    \FM^{(\lambda)}\}\}$.
\item{3.} If $\lambda$ is a limit ordinal then
          $\ds \FM^{(\lambda)}= \bigcap_{\mu < \lambda} \FM^{(\mu)}$.
\smallskip
We call the least $\lambda$ for which $\FM^{(\lambda)} \subseteq
\{ \emptyset \}$ the index of $\FM$ and denote it by $i(\FM)$.

\proclaim Definition 4. Let $I \subseteq \N$. Let $(\FM_k)_{k \in
I}$ be a sequence of compact hereditary families of finite subsets
of $\N$ and let $(\theta_k)_{k \in I} \subset (0,1]$. We denote by
$T \left[ \left( \FM_k , \theta_k \right)_{k \in I}\right] $ the
completion of $\c00$ with respect to the norm defined by $$ \norma
x =\max \left\{ \norma x_\infty, \sup_{k \in I} \left\{\theta_k
\sup_{n \in \N} \left\{ \sum_{i=1}^n \norma{E_ix}\vert\
(E_i)_{i=1}^n\ \FM_k-admissible \right\} \right\} \right\} $$ and
we call it the mixed Tsirelson space defined by the sequence
$(\FM_k,\theta_k)_{k \in I}$.

\proclaim Remark 2. The existence of such a norm is shown, for
instance, in [10]. It follows from the definition of the norm
that the sequence $(e_i)_{i=1}^\infty$ is a normalized
1-unconditional basis for $T \left[ \left( \FM_k , \theta_k
\right)_{k \in I}\right]$.

\proclaim Remark 3. There are two useful alternative ways to
define the norm. Given $\ds x= \sum_{n=1}^\infty a_n e_n \in
\c00$,
\item {(i)} define a non decreasing sequence of norms on $\c00$:
$$
\eqalign{
\vabs x_0 =& \max_{n \in \N} \vabs{a_n}
\cr
\vabs x_{s+1} =& \max
\left\{ \vabs x_s, \sup_{k \in I}
\left\{\theta_k \sup_{n \in \N}
\left\{ \sum_{i=1}^n \vabs{E_ix}_s \vert\ (E_i)_{i=1}^n\ \FM_k-admissible
\right\}
\right\}
\right\}
\cr}
$$
Then $\ds \norma x = \sup_{s \in \N \cup \{0\}} \vabs x_s$.
\item {(ii)}
Let $K_0= \{\pm e_n \ \vert\ n \in \N \}$. Given $K_s, s \in \N\cup \{0\}$, let
$$
\eqalign{
K_{s+1}= K_s \cup \Bigl\{ \theta_k& \cdot (f_1+\dots+f_d)\ \vert\
 k \in I,
d \in \N,
f_i \in K_s, i=1,\dots,d \cr &\hbox{ are successive and }
 (\sop(f_1),\dots,\sop(f_d))\ \FM_k-admissible
\Bigr\} \cr} $$ Let $\ds K=\bigcup_{s=0}^\infty K_s$. Then $\ds
\norma x = \sup \{f(x) \vert\ f \in K \}$.

The latter definition of the norm provides information about the
dual space. Looking at the set $K$ as a set of functionals it is
not difficult to see that $B_{X^*}$ is the closed convex hull of
$K$, where the closure is taken either in the weak-$*$ topology
or in the pointwise convergence topology.

\parrafo {1. The case $\theta_k=1$}

Let $J= \{k \in I \ \vert\ \theta_k =1 \}$. If $J$ is not empty,
we give information about the structure of $T \left[ \left( \FM_k
, \theta_k \right)_{k \in I}\right]$ depending on the index
$i(\FM_k), k\in J.$ It is known that if $i(\FM_k)\ge 2$ for some
$k\in J$, then $T \left[ \left( \FM_k , \theta_k \right)_{k \in
I}\right]$ contains an isomorphic copy of $\ell_1$. Actually it is
possible to say much more as our next proposition shows.

\proclaim Proposition 1. If $i(\FM_{k_0}) \geq 2$ for some
$k_0\in J$, then $T \left[ \left( \FM_k , \theta_k \right)_{k \in
I}\right] $ is $\ell_1-$saturated.

\dem By the Bessaga-Pelczynski principle (see e.g. [6], pg. 10),
it suffices to show that every block subspace contains a further
subspace isomorphic to $\ell_1$. Recall that a block subspace is
a space of the form $\overline {Span} \{ u_i, i\in\N\}$, with
$(u_i)_{i=1}^\infty$ a block sequence.

Let $(u_i)_{i=1}^\infty$ be a block sequence. We are going to
construct a subsequence $(u_{i_k})_{k=1}^\infty$ of
$(u_i)_{i=1}^\infty$ equivalent to the $\ell_1$ basis.

Let $\{p\} \in \FM_{k_0}^{(1)}$. We can choose $u_{i_1}$ such
that $p < u_{i_1}$. Now, since $\{p\} \in \FM_{k_0}^{(1)}$, there
exists $n_1 \in \N$ such that $n_1 > u_{i_1}$ and $\{p,n_1\} \in
\FM_{k_0}$, so we can take $u_{i_2}$ such that $n_1 < u_{i_2}$.
Continuing in this manner, we can construct a subsequence
$(u_{i_k})_{k=1}^\infty$ of $(u_i)_{i=1}^\infty$ such that for
every $k \in \N$ there exists $n_k \in \N$ such that $u_{i_k} <
n_k < u_{i_{k+1}}$ and $\{p, n_k\} \in \FM_{k_0}$. It is now easy
to see that $(u_{i_k})_{k=1}^\infty$ is equivalent to the
$\ell_1$ basis. \bull

 The following example shows a Tsirelson
type space $\ell_1$-saturated but not isomorphic to $\ell_1$. It
was shown to us by I. Deliyanni.

\proclaim Example 1. Let $\FM = \{ F \subseteq \N \ \vert \
\exists i \in \N \hbox { such that } F \subseteq \{1, 2^i\}\}$
and $\theta=1$.

 It is clear that $i(\FM)=2$. If $T[\FM,\theta]$
were isomorphic to $\ell_1$ then since $\ell_1$ has a unique --
up to equivalence -- normalized unconditional basis, there would
exist a constant $C> 0$ such that for all $n\in\N$, $$ {1 \over
C} \sum_{i=1}^n \vabs {a_i} \leq \norma {\sum_{i=1}^n a_i e_i}
\leq C  \sum_{i=1}^n \vabs {a_i}. $$

Now taking $\ds x = \sum_{i=2^k+1}^{2^{k+1}} e_i$ we would obtain
$2^k-1 \leq C$ for all $k \in \N$.
\medskip

We now examine $T \left[\left( \FM_k , \theta_k \right)_{k \in I}\right]$ with
$i(\FM_k)=1, k\in J$. We will find different subspaces depending on whether
 the set $\ds\bigcup_{k \in J} \FM_k$ contains only a finite number of non
singleton sets
 or not.\medskip

 \proclaim Proposition 2.  Let $I'\subseteq I$ be such
that $\ds\bigcup_{k \in I'} \FM_k$ contains only a finite number
of non singleton sets.
\item {(1)} If $I' \not= I$, then
$T \left[ \left( \FM_k , \theta_k \right)_{k \in I}\right]$ is isomorphic to
$T \left[ \left( \FM_k , \theta_k \right)_{k \in I\setminus I'}\right]$.
\item {(2)} If $I'= I$, then
$T \left[ \left( \FM_k , \theta_k \right)_{k \in I}\right]$ is isomorphic
to $c_0$.

\dem $(1).$ Let $\norma \cdot$ and $\norma \cdot'$ be the norms of the spaces
$T \left[ \left( \FM_k , \theta_k \right)_{k \in I}\right] $ and
$T \left[ \left( \FM_k , \theta_k \right)_{k \in I \setminus I'}\right] $,
respectively. We will see that they are equivalent. Clearly, $\norma \cdot'
\leq \norma \cdot$.

For the other inequality let  $\ds M=\max \left\{ \max E \ \vert \
E \in \bigcup_{k \in I'} \FM_k ,\, \hbox {non singleton} \right\}
$ and write \break $\ds x=\sum_{i=1}^\infty a_i e_i = \sum_{i=1}^M
a_i e_i + \sum_{i=M+1}^\infty a_i e_i := x_1 + x_2$.

We have $\norma {x_1} \leq M \norma x'$ since
$\ds
\norma {x_1}= \norma {\sum_{i=1}^M a_i e_i} \leq
 \sum_{i=1}^M \vabs {a_i} \leq
 \sum_{i=1}^M \norma x_\infty \leq M \norma x'.$

On the other hand, we show first by induction over $s$ that $
\vabs {x_2}_s \leq \vabs {x_2}_s'$. For $s=0$ it is clear. Suppose
now that it is true for $s$ and let
 $E_1, \dots, E_n$ be a sequence of finite subsets of $\N$,
$\FM_k-$ admissible for some $k$. There are two possibilities,
either  $k \in I \setminus I'$ and then
$\ds \theta_k \sum_{i=1}^n \vabs {E_ix_2}_s \leq
     \theta_k \sum_{i=1}^n \vabs {E_ix_2}_s' \leq \vabs
     {x_2}_{s+1}'$, or $k \in I'$ and then, by hypothesis, $n=1$, $E_1$ is
$\FM_k-$admissible and $ \theta_k \vabs {E_1x_2}_s \leq \theta_k
\vabs {x_2}_s \leq \vabs {x_2}_s' \leq \vabs
{x_2}_{s+1}'$.\smallskip

Therefore, $\norma {x_2} \leq \norma{x_2}'$ and by
$1-$unconditionality, $\norma {x_2}' \leq \norma x'$. Thus, $
\norma x' \leq \norma x \leq (M+1) \norma x'$.

For $(2)$, it is easy to see that $T \left( \FM_0 , \theta_0
\right)$ is isomorphic to $c_0$, where $\FM_0 = \{ \{i\} \ \vert\
i \in \N \}$, and $\theta_0=1$. Now use (1) to  get that $ T
\left[ \left( \FM_k , \theta_k \right)_{k \in I}\right]$ is
isomorphic to $T \left[ \left( \FM_k , \theta_k \right)_{k \in I
\cup\{0\} }\right]$ and once again to see that the latter is
isomorphic to $ T \left( \FM_0 , \theta_0 \right).$ \bull

Proposition 2 for $I'=J$ yields

\proclaim Proposition 3. Let $J= \{k \in I \ \vert\ \theta_k=1\}$.
\item {(1)} Let $\ds \bigcup_{k \in J} \FM_k$ contain only a finite number
of non singleton sets.
\itemitem {1.1.} If $J = I$, then $T \left[\left(\FM_k, \theta_k\right)_{k\in
I}\right]$ is isomorphic to $c_0$.
\itemitem {1.2.} If $J \not= I$, then
$T \left[ \left( \FM_k , \theta_k \right)_{k \in I}\right]$ is isomorphic to
$T \left[ \left( \FM_k , \theta_k \right)_{k \in I\setminus J}\right]$.
\item {(2)} Let  $\ds \bigcup_{k \in J} \FM_k$ contain an infinite number of
non singleton sets. Then $T \left[ \left( \FM_k , \theta_k
\right)_{k \in I}\right]$ contains  a subspace isomorphic to
$\ell_1$.

\dem $(1)$ follows from Proposition 2. For $(2)$, we will
construct a subsequence $(e_{n_i})_{i=1}^\infty$ of
$(e_i)_{i=1}^\infty$ equivalent to the $\ell_1$ basis.

Let $\ds M_1 \in \bigcup_{k \in J} \FM_k$ be a non singleton. Let
$n_1=\min M_1$. Having chosen $n_i$, we can take $\ds M_{i+1} \in
\bigcup_{k \in J} \FM_k$ a non singleton such that $\min M_{i+1}
> \max M_i$, and take $n_{i+1}=\min M_{i+1}$.

Consider the sequence $(e_{n_i})_{i=1}^\infty$ and let's show
that it is equivalent to the $\ell_1$ basis.

Let $\ds x=\sum_{i=1}^\infty a_i e_{n_i}$. By the definition of
the norm and the fact that for every $N \in \N$ and $i < N$,
$(\{n_i\}, [n_{i+1}, n_N]\cap \N \})$ is $\FM_k-$admissible for
some $k \in J$ we have $$ \norma x \geq \vabs {a_1} + \norma
{\sum_{i=2}^N a_i e_{n_i}} \geq\dots\geq
 \vabs {a_1}+\vabs {a_2} + \dots + \vabs {a_N}.$$

The proof is complete since always $\norma x \leq \norma x_1$.
\bull

Observe that in statement (2) of Proposition 3 we do not ensure
$\ell_1$ saturation. Actually, in some cases we can also find
$c_0$ as a subspace. This is a consequence of the following
general result.

\proclaim Proposition 4. Let $\FM_k$ be compact and hereditary
for all $k \in I \subseteq \N$, $\theta_k \in (0,1]$ for all $k
\in I$. If for all $N \in \N$ there exists $n \geq N$ such that
for all $\ds M \in \bigcup_{k \in I}\FM_k$ either $\ n < \min M$
or $n \geq \max M$, then $T \left[ \left( \FM_k , \theta_k
\right)_{k \in I}\right] $ contains a subspace isomorphic to
$c_0$. Moreover, if $\theta_k=1$ for all $k \in I$, the converse
is true.

\dem We will construct a subsequence $(e_{n_i})_{i=1}^\infty$ of
the basis $(e_i)_{i=1}^\infty$ equivalent to the basis of $c_0$.

Let $N_1=1$. By hypothesis there exists $n_1 \geq N_1$ such that
for all $\ds M \in \bigcup_{k \in I} \FM_k$, $n_1 < \min M$ or
$n_1 \geq \max M$.

Suppose that $n_i$ is chosen and write $N_{i+1}=n_i+1$. Then
there exists $n_{i+1} \geq  N_{i+1}$ verifying the hypothesis.
Now, consider the sequence $(e_{n_i})_{i=1}^\infty$.

Let $\ds x= \sum_{i=1}^\infty a_i e_{n_i}\in c_{00}$
 and write $\vabs x_0 = \norma x_\infty$ as in Remark 3.

Let $(E_i)_{i=1}^n$ be a sequence of finite subsets of $\N$, $\FM_k-$
admissible for some $k \in I$. Then we have $\ds \theta_k \sum_{i=1}^n
\vabs {E_ix}_0 = \theta_k \vabs {E_{i_0}x}_0 \leq \vabs x_0$ and so
$\vabs x_1 \leq \vabs x_0$. Indeed, the first equality is true
since by the construction of $(n_i)$, there exists at most one $E_i$ such
that $\sop (x) \cap E_i \not= \emptyset$ and the inequality is
straightforward by 1-unconditionality. So we have proved that
$\vabs x_1 = \vabs x_0$ and therefore $\vabs x_n = \vabs x_{n+1}$
and $\norma x = \norma x_\infty$.
\medskip

The converse is a consequence of the following
\medskip

{\it CLAIM:} {\sl If there is an $N_0$ such that for all $n \geq
N_0$, there exists $\ds M \in \bigcup_{k \in I} \FM_k$ such that
$\min M \leq n < \max M$, then every normalized block sequence in
$T[(\FM_k,1)_{k\in I}]$ has a subsequence equivalent to the
canonical basis of $\ell_1$ and in particular, $T[(\FM_k,1)_{k\in
I}]$ is $\ell_1-$saturated.}
\medskip

{\it Proof of CLAIM:}
 Let $(x_i)_{i=1}^\infty$ be a normalized block sequence. Let $i_1$ be such
 that $N_0 \leq \min x_{i_1}$.
We split $\ds x_{i_1}= \sum_{k=p_1+1}^{p_2} a_k e_k$ in the
following manner:

Let $\ds A^{(1)}(x_{i_1}) = \left\{ j > \min x_{i_1} \ \vert\
\{t,j\} \in \bigcup_{k \in I} \FM_k\,,t \leq \min x_{i_1}
\right\}$. By hypothesis  $A^{(1)}(x_{i_1})$ is not empty and
$j^{(1)}(x_{i_1}) := \min A^{(1)}(x_{i_1}) > \min x_{i_1}$.

Therefore, $$ x_{i_1}= \sum_{k=p_1+1}^{p_2} a_k e_k =
\sum_{k=p_1+1}^{j^{(1)}(x_{i_1})-1} a_k e_k
+\sum_{k=j^{(1)}(x_{i_1})}^{p_2} a_k e_k := x_{i_1}^{(1)} +
u^{(1)}. $$

Let $\ds y_{i_1}^{(1)}= {x_{i_1}^{(1)} \over \norma
{x_{i_1}^{(1)}}}$. Suppose $y_{i_1}^{(\ell)}$ is defined and we
have $x_{i_1}= x_{i_1}^{(1)} + \dots + x_{i_1}^{(\ell)} +
u^{(\ell)}$. If $u^{(\ell)} \not= 0$, define $x_{i_1}^{(\ell+1)} =
\left(u^{(\ell)} \right)^{(1)}$ and $\ds y_{i_1}^{(\ell+1)}=
{x_{i_1}^{(\ell+1)} \over \norma {x_{i_1}^{(\ell+1)}}}$ and keep
going until we have $u^{(d_1)}=0$ for some $d_1\in\N$. Then we
have $\ds x_{i_1} = \sum_{\ell=1}^{d_1} \norma {x_{i_1}^{(\ell)}}
y_{i_1}^{(\ell)}. $

Now, take $i_2$ such that $\sop (x_{i_2}) > j^{(d_1)}(x_{i_1})$
 and split it as before. Continuing in this manner, we obtain a sequence
$$
\left( y_{i_1}^{(1)}, y_{i_1}^{(2)}, \dots , y_{i_1}^{(d_1)}, y_{i_2}^{(1)},
 \dots, y_{i_2}^{(d_2)}, \dots ,y_{i_n}^{(1)}, \dots, y_{i_n}^{(d_n)}, \dots
\right) := (u_k)_{k=1}^\infty. $$

For this sequence we have $$ \norma{\sum_{k=1}^n a_ku_k}=
\vabs{a_1}+\norma{\sum_{i=2}^n a_ku_k}= \dots
=\sum_{k=1}^n\vabs{a_k}, $$ that is, $(u_k)_{k=1}^\infty$ is
equivalent to the canonical basis of $\ell_1$. But
$(x_{i_k})_{k=1}^\infty$ is a block sequence of
$(u_k)_{k=1}^\infty$ and therefore it is also equivalent to the
canonical basis of $\ell_1$. \bull

\proclaim Remark 4.
\item {1.} Observe that, in particular, the hypothesis of
Proposition 4 implies that $i(\FM_k)=1$ for all $k \in I$.
\item {2.} The proof of the converse of Proposition 4 states that either
 $T[(\FM_k,1)_{k\in I}]$ contains a subspace
isomorphic to $c_0$ or
 $T[(\FM_k,1)_{k\in I}]$ is $\ell_1-$saturated.
 \smallskip

We now give an example of a Tsirelson type space which contains $\ell_1$
and $c_0$.

\proclaim Example 2. Let $\FM = \{ F
\subseteq \N \ \vert \ \exists i \in \N \hbox { such that } F \subseteq
\{2i-1, 2i\}
 \}$. $T(\FM,1)$ contains
$\ell_1$ by Proposition 3 and $c_0$ by Proposition 4. Moreover, it
is easy to see that the space is isomorphic to $\ell_1 \oplus
c_0$.

\parrafo {2. The case $(\FM_k,\theta_k)_{k=1}^\ell$}
In view of the previous results, in this section we will consider
Tsirelson type spaces defined by finite sequences
$(\FM_k,\theta_k)_{k=1}^\ell$,  with $\theta_k \in (0,1)$ for all
$k=1,\dots,\ell$. The main result of the section is

\proclaim Theorem 1. Let $i(\FM_k) = n_k\in\N$ and $\theta_k \in
(0,1)$ for all $k=1,\dots,\ell$.
\item {1.} If $\ds \theta_k \leq {1 \over n_k}$ for all $k$ then
$T[(\FM_k,\theta_k)_{k=1}^\ell]$ is
$c_0-$saturated.
\item {2.} If $\ds \theta_k > {1 \over n_k}$ for some $k$ then
$T[(\FM_k,\theta_k)_{k=1}^\ell]$ is $\ell_p-$saturated for some
$p\in(1,+\infty)$.

Our proof of this theorem is based on Theorem 2 below, proved in
[4]. In order to state it we first need some definitions.

\proclaim Definition 5. Let $m\in\N$ and $\phi\in K_m\setminus
K_{m-1}$. An analysis of $\phi$ is any sequence
$\{K^s(\phi)\}_{s=0}^m$ of subsets of  $K$ such that for every
$s$,
\item  {1.} $K^s(\phi)$ consists of successive elements of  $K_s$ and
$\ds \bigcup_{f\in K^s(\phi)} \sop(f)= \sop(\phi)$.
\item  {2.} If $f\in K^{s+1}(\phi)$ then either $f\in K^s(\phi)$ or
there exists $k$ and successive  $f_1,\dots,f_d \in K^s(\phi)$
with $(\sop(f_1),\dots, \sop(f_d))$ $\FM_k-{\rm admissible}$ and
$f=\theta_k(f_1+\dots+f_d)$.
\item  {3.} $K^m(\phi)=\{\phi\}$.\eject

 \proclaim Definition 6.
\item {1.} Let $\phi\in K_m\setminus K_{m-1}$ and let
$\{K^s(\phi)\}_{s=0}^m$ be a fixed analysis of  $\phi$. Then for a
given finite block sequence $(x_k )_{k=1} ^{\ell}$ we set for
every $k \in \{1,\dots,\ell\}$ $$ s_k =\cases {
 \max \{ s\ \vert\ 0\leq s <m \hbox {, and there are  at least two } f_1,
f_2 \in  K^s(\phi) \cr \phantom {\max \{}\hbox { such that } \vabs
{\sop(f_i)\cap \sop(x_k)} >0,
         i=1,2\},\cr
\phantom {\max}  \hbox {when this set is non - empty}&\cr \quad\cr
0 \qquad \hbox{if } \vabs {\sop (x_k)\cap \sop(\phi)}\leq 1. &\cr}
$$
\item {2.} {For $k=1,\dots , \ell $ we define the initial part and
the final part of $x_k$ with respect to $\{ K^s(\phi )\}_{s=0} ^m
$, and denote them respectively by $x'_k\, $ and $x''_k,\,$ as
follows: If $\{f\in K^{s_k}(\phi)\vert \sop(f)\cap
\sop(x_k)\ne\emptyset\}:=\{f_1,\dots ,f_d\}$ with $f_1<\dots<f_d$,
we set $x'_k = (\sop(f_1))x_k$ and $\ds x''_k =\left(\cup _{i=2}^d
\sop(f_i)\right)x_k$. }

\proclaim Notation. Let  $m \in \N,\ \ \phi \in K^m \setminus
K^{m-1}$, let $\{K^s(\phi)\}_{s=0}^m$ be an analysis of $\phi$,
$(v_i)_{i=1}^\infty$ a block sequence and $(x_j)_{j=1}^\infty$ a
block sequence with $x_j\in Span\{v_i\mid i\in\N\}$. Suppose that
there exists $n_\phi$ such that $\ds \sop (\phi) \subseteq
\bigcup_{j=1}^{n_\phi} \sop (x_j)$ and denote by $x_j'$ and
$x_j''$ the initial and the final part of $x_j$, $j \leq n_\phi$.
For all $f = \theta_k(f_1+\dots+f_d) \in K^s(\phi)$ and $J
\subseteq \{1,\dots,n_\phi\}$ we define the following sets for
$(x_j')$: $$ I'= \{i \ \vert\ 1 \leq i \leq d \hbox{ and } \sop
(f_i) \cap \sop (x_j') \not= \emptyset \hbox{ for at least two
different } j \in J\} $$ and for every $i \in I$, $$D'_{f_i}= \{
j \in J \ \vert \ \sop (f_i) \cap \sop (x_j') \not= \emptyset
\hbox { and } (\sop (f) \cap \sop (x_j') ) \setminus \sop (f_i)
\subseteq \sop (v_t) \hbox { for some $t$} \} $$ and $$
\eqalign{T'= \{j \in J\ & \vert \ j \notin \bigcup_{i \in I'}
D'_{f_i}\hbox { and } \exists t_1 \not= t_2 \cr &\hbox{ such that
} \sop(x_j') \cap \left( \cup_{i \notin I'} \sop (f_i) \right)
\cap \sop (v_{t_i}) \not= \emptyset, \ \ \  i=1,2 \}. \cr} $$

In the same manner we define sets $I'',D''_{f_i}, T''$ exchanging
$x_j'$ for $x_j''$.

\proclaim Theorem 2 ([4]).  Given
$T[(\FM_k,\theta_k)_{k=1}^\ell]$ with $\ell \in \N$, $\theta_k \in
(0,1)$ and $i(\FM_k) = n_k\in\N$, for all $k=1,\dots,\ell$, let
$(v_i)_{i=1}^\infty$ be a normalized block sequence. If there
exists a sequence $\ds x_j=\sum_{i \in I_j} a_i v_i$ with
$(a_i)_{i=1}^\infty\subset\R$ and $(I_j)_{j=1}^\infty\subset\N$
successive such that
\item {(a)} $\ds {1 \over 2^{j+1}}\leq \vabs {a_j} < {1  \over 2^j}$
and
\item {(b)} for all $m \in \N,\ \ \phi \in K^m \setminus K^{m-1}$,
each analysis $\{K^s(\phi)\}_{s=1}^m$ of $\phi$, all $f=\theta_k
(f_1+\dots+f_d) \in K^s(\phi)$, and all $J \subseteq
\{1,\dots,n_\phi\}$, the inequalities $\vabs{I'} + \vabs{T'} \leq
n_k$ and  $\vabs{I''} + \vabs{T''} \leq n_k$ hold,\smallskip
\noindent then $(x_j)_{j=1}^\infty$ is equivalent to the canonical
basis of $T[(\FA_{n_k},\theta_k)_{k=1}^\ell]$.
\medskip\medskip

Recall, see [4], that the space
$T[(\FA_{n_k},\theta_k)_{k=1}^\ell]$ is either isometrically
isomorphic to $c_0$, when $n_k \cdot \theta_k \leq 1$ for all $k$,
or isomorphic to $\ell_p$, where $\ds p=min \left\{ {1 \over
1-\log_{n_k} {1 \over \theta_k}} \ \vert \ n_k \cdot \theta_k >1
\right\}.$ So, to prove Theorem 1 we need to find the sequence
$(x_j)_{j=1}^\infty$  and the next lemma will be useful for
constructing it.
\smallskip

\proclaim Lemma 1. Let $\ell \in \N$, $\theta_k\in(0,1)$ and
$\FM_k$ be such that $i(\FM_k)=n_k\in\N$ for all $k=1,\dots\ell$.
Then for every  block sequence $(u_i)_{i=1}^\infty$ in
$T[(\FM_k,\theta_k)_{k=1}^\ell]$ there exists an infinite subset
$\FP=\{p_i\}_{i=1}^\infty$ of $\N$ and a subsequence
$(v_i)_{i=1}^\infty$ of $(u_i)_{i=1}^\infty$ having the following
properties:
\item {(a)}$ p_1 \leq \sop (v_1) < p_2 \leq \sop (v_2) < \dots < p_i \leq
\sop (v_i) < p_{i+1} \leq
\dots$
\item {(b)} For every sequence
$E_1 <E_2 \dots < E_{n_k}$ of finite subsets of $\FP$, where $E_i
= \{ p_{\ell_1^i}, \dots , p_{\ell_{t_i}^i}\}, i=1,\dots, n_k$,
the family $$ \left( \bigcup_{j=\ell_1^1}^{\ell_{t_1}^1} \sop
(v_j), \dots,
       \bigcup_{j=\ell_1^{n_k}}^{\ell_{t_{n_k}}^{n_k}} \sop (v_j) \right)
$$
is $\FM_k-$admissible.
\item {(c)} If $r \geq n_k+1,\, S=\{s_1, \dots s_r\} \subseteq \N$
is such that $$
 \vabs { \{j \in \N\ \vert\ [s_i, s_{i+1}] \cap \sop (v_j) \not=
\emptyset \} } \geq 2 $$ for all $i=1, \dots, r-1$, then $S \notin
\FM_k$.

\dem  The proof is based on the following result from
[4]:\smallskip

\proclaim Lemma 2. Let $\ell,n_1,\dots,n_\ell \in \N$. Let
$\FM_k,k=1,\dots, \ell$ be such that $i(\FM_k)=n_k$. Then there
exists an infinite subset $Q$ of $\N$ having the following
properties:
\item {1.} Let $k \in \{1,\dots,\ell\}$. Every sequence
$E_1 <E_2 \dots < E_{n_k}$ of length $n_k$ of finite subsets of
$Q$ is $\FM_k-$admissible.
\item {2.} Let $k \in \{1,\dots,\ell\}$. If $r \geq n_k +1$, then
no sequence $E_1 < E_2 \dots < E_r$ of finite subsets of $Q$ with
$\vabs {E_i} \geq 2$ for all $i=1,\dots,r$, is
$\FM_K-$admissible.\smallskip\smallskip

Now, let $Q=\{k_i \}_{i=1}^\infty$ be the sequence in Lemma 2.
Take $p_1=k_1$, and $v_1=u_{\ell}$ such that $p_1 \leq \sop
(u_{\ell})$. Having chosen $p_i$ and $v_i$ with $p_i \leq \sop
(v_i)$, since $\{k_i \}_{i=1}^\infty$ is increasing, let
$k_{j_i}$ be such that $p_i \leq \sop (v_i) < k_{j_i}$, and take
$p_{i+1} = k_{j_i+1}$ and $v_{i+1} = u_{\ell}$ such that $p_{i+1}
\leq \sop (u_{\ell})$.

The sequences $\{p_i \}_{i=1}^\infty$ and $(v_i)_{i=1}^\infty$
satisfy the assertions of Lemma 1:
\item {(a)} By construction.
\item {(b)} It is sufficient to see that
$\ds \bigcup_{j=\ell_1^i}^{\ell_{t_i}^i} \sop (v_j) \subseteq
\left[ p_{\ell_1^i}, p_{\ell_{t_i}^i}\right]$ and, since the
family $\ds \left\{ \left\{p_{\ell_1^i}, p_{\ell_{t_i}^i} \right\}
\right \}_{i=1}^{n_k}$ is $\FM_k-$ad\-mi\-ssi\-ble by Lemma 2,
$\ds \left(\bigcup_{j=\ell_1^i}^{\ell_{t_i}^i} \sop (v_j)
\right)_{i=1}^{n_k}$ is also $\FM_k-$admissible.
\item {(c)} Let $S=\{ s_1, \dots ,s_r\}$ be such that $
\vabs { \{j \in \N\ \vert\ [s_i, s_{i+1}] \cap \sop (v_j) \not=
\emptyset \} } \geq 2$ for all $i=1,\dots,r-1$,  let $d_i = \min
\{j \in \N \ \vert \ [s_i, s_{i+1}] \cap \sop (v_j) \not=
\emptyset \}$. Then $k_{j_{d_i}}$ and $ p_{d_i+1} \in [s_i,
s_{i+1}] \cap Q$ for all $i=1, \dots, r-1$, and by the property
(2) of Lemma 2, $S \notin \FM_k$. \bull
\bigskip\medskip

\dem {\sl (of Theorem 1).} It suffices to show that $c_0$ or
$\ell_p$ is included in every block subspace of
$T[(\FM_k,\theta_k)_{k=1}^\ell]$.

Let  $(u_i)_{i=1}^\infty$ be a normalized block sequence. Let
$\FP=\{p_i\}_{i=1}^\infty$ and $(v_i)_{i=1}^\infty$ be the
sequences associated to $(u_i)_{i=1}^\infty$ from Lemma 1.

 If
$\ds \sup_{m \in \N} \norma{\sum_{i=1}^m v_i}$ is finite, then
$(v_i)_{i=1}^\infty$ is equivalent to the canonical basis of
$c_0$, and from Corollary 1 of [4] we have $n_k \cdot \theta_k
\leq 1$.

Suppose now that  $\ds \lim_{m \to \infty} \norma{\sum_{i=1}^m
v_i}=\infty$. Then we can construct a sequence
$(y_i)_{i=1}^\infty$ supported by the subsequence
$(v_i)_{i=1}^\infty$ with the following properties: For every
$\ds j,y_j= {1 \over 2^{j+1}} \sum_{i \in I_j} v_i$, where
\item {(i)} $I_j$ are successive intervals of $\N$, and
\item {(ii)} $\ds 1-{1 \over {2^{j+1}}} \le \norma {y_j}\leq 1$.

If $x_j=\ds { y_j \over \norma {y_j}}$, the sequence $x_j$
satisfies condition (a) of Theorem 2.

We prove condition (b) of Theorem 2 for the initial parts of
$(x_j)$ since for the final parts the proof is analogous. Suppose
that $\phi, f$ and $J$ are fixed. Let $m_1 \leq \sop (f_1) < m_2
\leq \sop(f_2) < \dots<m_d\leq \sop (f_d)$. We define $B
\subseteq \{m_1,\dots,m_d \}$ as follows:
$$ m_{i_s} \in B \ssi
\cases { i_s \in I', &\cr i_s= \min \{i \notin I' \ \vert \
\sop(x_j') \cap \sop (f_i) \not= \emptyset \} \hbox { for some }
j \in T'. & \cr} $$

Let $m_{i_1} < \dots < m_{i_r}$ be the elements of $B$. Observing
that $$ \vabs { \left\{ t \in \N \ \vert \ [m_{i_s}, m_{i_{s+1}}]
\cap \sop (v_t) \right\}} \geq 2,\qquad \forall\ 1 \leq s\leq
r-1$$
and using property (c) of Lemma 1 we get that $r= \vabs B
\leq n_k$. So $\vabs {I'} + \vabs {T'} \leq n_k$. \bull

The proof of the next two corollaries easily follows from
 Theorem 1 from this paper and
Corollaries 1 and 2 from [4].

\proclaim Corollary 1. Let $T[(\FM_k,\theta_k)_{k=1}^\ell]$,
$1<p<\infty$, $n_k=i(\FM_k)$ and $\theta_k \in (0,1)$. The
following conditions are equivalent:
\item {i)} $T[(\FM_k,\theta_k)_{k=1}^\ell]$ contains a
subspace isomorphic to $\ell_p$.
\item {ii)} $T[(\FM_k,\theta_k)_{k=1}^\ell]$ is $\ell_p-$saturated.
\item {iii)} $i(\FM_k)$ is finite, $\ds \theta_k > {1 \over n_k}$ for
some $k =1,\dots,\ell$ and $\ds p= \min \left\{{1 \over
1-\log_{n_k} {1 \over \theta_k}} \ \vert\ n_k \cdot \theta_k > 1
\right\}$.

\proclaim Corollary 2. Let $T[(\FM_k,\theta_k)_{k=1}^\ell]$,
$\theta_k \in (0,1)$. The following conditions are equivalent:
\item {i)} $T[(\FM_k,\theta_k)_{k=1}^\ell]$ contains a
subspace isomorphic to $c_0$.
\item {ii)} $T[(\FM_k,\theta_k)_{k=1}^\ell]$ is $c_0-$saturated.
\item {iii)} $i(\FM_k)$ is finite and $\ds \theta_k \leq {1 \over i(\FM_k)}$
for all $k =1,\dots,\ell$.

In view of Proposition 1 and the previous corollaries we can
include the case $\ell_1$ in the discussion.

 \proclaim Corollary 3. Let
$T[(\FM_k,\theta_k)_{k=1}^\ell]$, $2\le i(\FM_k)\in\N$  and
$\theta_k\in (0,1]$.  The following conditions are equivalent:
\item {i)} $T[(\FM_k,\theta_k)_{k=1}^\ell]$ contains a
subspace isomorphic to $\ell_1$.
\item {ii)} $T[(\FM_k,\theta_k)_{k=1}^\ell]$ is $\ell_1-$saturated.
\item {iii)} $\theta_k =1$ for some $k =1,\dots,\ell$

So in particular we have proved the following criterion, which is
useful to show when two Tsirelson type Banach spaces are totally
incomparable.

\proclaim Theorem 3. Let $\ell, \ell' \in \N$, $\theta_k \in
(0,1)$ and $i(\FM_k)=n_k\in\N$ for all $k = 1, \dots,\ell$ and
$\theta'_k \in (0,1)$ and $i(\FM'_k)=n'_k\in\N$ for all $k
=1,\dots,\ell'$. Then $T[(\FM_k,\theta_k)_{k=1}^\ell]$ and
$T[(\FM'_k,\theta'_k)_{k=1}^{\ell'}]$ are totally incomparable if
and only if one of the following situations occurs:
\item {1.} $\ds \theta_k \leq {1 \over n_k}$ for all $k=1,\dots,\ell$ and $\ds
\theta'_k > {1 \over n'_k}$ for some $k \in \{1,\dots,\ell'\}$, or
\item {2.} $\ds \theta'_k \leq {1 \over n'_k}$ for all $k=1,\dots,\ell'$ and
$\ds \theta_k > {1 \over n_k}$ for some $k \in \{1,\dots,\ell\}$,
or
\item {3.} $\ds \theta_k > {1 \over n_k}$ for some $k \in
\{1,\dots,\ell\}$ and $\ds \theta'_k > {1 \over n'_k}$ for some $k
\in \{1,\dots,\ell'\}$ and $$
 \min \left\{{1 \over 1-\log_{n_k} {1 \over \theta_k}} \ \vert\
n_k \cdot \theta_k > 1 \right\} \not= \min \left\{{1 \over
1-\log_{n'_k} {1 \over \theta'_k}} \ \vert\ n'_k \cdot \theta'_k
> 1\right\}.$$

Also we obtain a characterization of the reflexivity of this kind
of spaces as in [1].

\proclaim Proposition 5. Let $\ell \in \N$. Let $\theta_k \in
(0,1)$ and $i(\FM_k) = n_k\in\N$ for all $ k=1,\dots \ell$. Then
the following conditions are equivalent:
\item {1.} $T[(\FM_k, \theta_k)_{k=1}^\ell]$ is reflexive.
\item {2.} $\ds \theta_k > {1 \over i(\FM_k)}$ for some $k \in \{1, \dots
,\ell\}$.

\parrafo {3. A criterion of total incomparability for spaces of the form
$T[(\FA_k,\theta_k)_{k=1}^\infty]$} We will suppose throughout the
section that $(\theta_k)_{k=1}^\infty\subset(0,1] $ is a non
increasing null sequence since $T[(\FA_k,\theta_k)_{k=1}^\infty]$
is easily seen to be isometric to
$T[(\FA_k,\theta'_k)_{k=1}^\infty]$ where
$\theta_k'=\sup\{\theta_j\mid j\ge k\}$ and $\inf\{\theta_k\}>0$
implies that $T[(\FA_k,\theta_k)_{k=1}^\infty]$ is isomorphic to
$\ell_1$.

The following properties of such spaces, stated as lemmas, are
known.

\proclaim  Lemma 3. Let $(u_i)_{i=1}^n$ be a normalized block
sequence in $T[(\FA_k,\theta_k)_{k=1}^\infty]$. Then for all
$a_1,\dots,a_n \in \R$, $$ \norma {\sum_{i=1}^n a_i e_i} \leq
\norma {\sum_{i=i}^n a_iu_i}. $$

\dem It is easy to prove by induction on $s$ that $\ds
\vabs{\sum_{i=1}^n a_ie_i}_s \leq \norma {\sum_{i=i}^n a_iu_i}$.
\bull

The following lemma was proved in [11] with $\theta_k =
(\log_2(1+k))^{-1}$, but the same proof works for any $\theta_k$
converging to zero.

\proclaim Lemma 4 ([11]). Let
$T\left[(\FA_k,\theta_k)_{k=1}^\infty\right]$, let $\theta_k$
converge to $0$. Let $(y_n)_{n=1}^\infty$ be a block sequence,
let a strictly decreasing null sequence
$(\varepsilon_n)_{n=1}^\infty \subset \R^+$ and a strictly
increasing sequence $(k_n)_{n=1}^\infty \subset \N$ be such that
for each
 $n$ there is a normalized block sequence $(y(n,i))_{i=1}^{k_n},\
(1+\varepsilon_n)-$equivalent to the $\ell_1^{k_n}$ basis and $\ds
y_n={1 \over k_n} \sum_{i=1}^{k_n} y(n,i)$. Then for all $\ell \in
\N$, $$ \lim_{n_1 \to \infty} \lim_{n_2 \to \infty} \dots
\lim_{n_\ell \to \infty} \norma {\sum_{i=1}^\ell y_{n_i}} = \norma
{\sum_{i=1}^\ell e_i}. $$

We will consider spaces such that $\ell_1$ is finitely block
represented in every block subspace of the space but not
containing $\ell_1$. The role of $\ell_1$ in this context, as
well as that of $c_0$, can be easily described:

\proclaim Proposition 6. The following conditions are equivalent:
\item{i)} The identity is an isometric isomorphism from
$T[(\FA_k,\theta_k)_{k=1}^\infty]$ onto $c_0$.
\item{ii)} $T[(\FA_k,\theta_k)_{k=1}^\infty]$
 contains a subspace isomorphic to $c_0$.
\item{iii)} For all $n\in\N$, $\Vert\sum_{i=1}^{n} e_i\Vert=1$.
\item{iv)} $\theta_k\le 1/k$ for all $k\in\N$.

\dem $ii)\Rightarrow iii)$: By the Bessaga-Pelcynski Principle and
a theorem of R.C. James (see e.g. [8], pg. 97), for every
$\varepsilon>0$ there exists a normalized block sequence
$(u_i)_{i=1}^\infty$ such that for all $\ell \in\N$,
$$\max|a_i|\le \Vert\sum_{i=1}^{\ell} a_i u_i\Vert\le
(1+\varepsilon)\max|a_i|\qquad a_1\dots a_{\ell}\in\R$$ and so by
Lemma 3, $\ds \Vert\sum_{i=1}^{\ell} e_i\Vert\le (1+\varepsilon)$
and iii) follows. $iii)\Rightarrow iv)$: This is clear since
$\ds\theta\cdot\ell\le \Vert\sum_{i=1}^{\ell} e_i\Vert$.
$iv)\Rightarrow i)$: By induction on $m\in\N$ it easily follows
that $|\cdot|_m=|\cdot|_0$ on $c_{00}$. \bull

\proclaim Proposition 7. Let
$T\left[(\FA_k,\theta_k)_{k=1}^\infty\right]$, let $\theta_k$
converge to $0$ . The following conditions are equivalent:
\item{i)} The identity is an isometric isomorphism from
$T\left[(\FA_k,\theta_k)_{k=1}^\infty\right]$ onto $\ell_1$.
\item{ii)} $T\left[(\FA_k,\theta_k)_{k=1}^\infty\right]$
 contains a subspace isomorphic to $\ell_1$.
\item{iii)} For all $n\in\N$, $\Vert\sum_{i=1}^{n} e_i\Vert=n$.
\item{iv)} $\theta_2=1$.

\dem $ii)\Rightarrow iii)$. Choose a strictly decreasing sequence
$(\varepsilon_n))_{n=1}^\infty\subset\R_{+}$ converging to $0$ and
$k_n=n$. We will construct a block sequence $(y_n)_{n=1}^\infty$
as in Lemma 4 above.

By James' Theorem let $(u_i))_{i=1}^\infty$ be a normalized block
sequence $(1+\varepsilon_1)-$equivalent to the unit vector basis
of $\ell_1$. Let $y_1=u_1$. Again by James' theorem there exist a
normalized block sequence $(u_i'))_{i=1}^\infty$  with $u_i'\in
Span\{u_i\mid i\in\N\}$ and $y_1<u_1'$,
$(1+\varepsilon_2)-$equivalent to the unit vector basis of
$\ell_1$. Let $y_2={1\over 2}(u_1'+u_2')$. We continue in the
same way.

Let $\ell\in\N$. Since any block sequence
$(y_{n_i})_{i=1}^{\ell}$ is $(1+\varepsilon_1)-$equivalent to the
unit vector basis of $\ell_1^{\ell}$, by Lemma 4 we have
$$(1-\varepsilon_1)\ell\le\Vert\sum_{i=1}^{\ell}e_i\Vert\le \ell$$
and the result follows. $iii)\Rightarrow iv):$ Just notice that
$2=\Vert e_1+e_2\Vert=2\theta_2$. $iv)\Rightarrow i):$ This
follows by induction on $|\sop(x)|$. \bull

We now give sufficient conditions, in terms of the behaviour of
$\ds \lambda_n :=\norma{\sum_{i=1}^n e_i}$, guaranteeing that in
a space of this kind $\ell_1$ is finitely block represented in
every block subspace.

\proclaim Proposition 8 ([5]).  Let $n,\ell \in\N,
0<\varepsilon<1$. Let $(X,\norma \cdot )$ be a normed space with
a normalized 1-unconditional normalized basis
$(e_i)_{i=1}^{n^\ell}$ such that $$ (n-\varepsilon)^\ell\leq
\norma { \sum_{i=1}^{n^\ell}e_i}\leq n^\ell. $$ Then there exists
a normalized block sequence $(y_i)_{i=1}^n$ of
$(e_i)_{i=1}^{n^\ell}$ such that $$ n-\varepsilon \leq \norma
{\sum_{i=1}^n y_i}\leq n. $$ Moreover, $(y_i)_{i=1}^n$ is ${1\over
1-\varepsilon}$-equivalent to the canonical basis of $\ell_1^n$.

\proclaim Proposition 9. Let $T[(\FA_k,\theta_k)_{k=1}^\infty]$,
let $\theta_k$ converge to $0$. If there exists
$(n_k)_{k=1}^\infty \subseteq \N$ unbounded and
$(\ell_k)_{k=1}^\infty$ such that $$ \lim_{k\to \infty} \left[n_k
- \left(\lambda_{n_k^{\ell_k}}\right)^{1 \over \ell_k} \right] =
0, $$ then $\ell_1$ is finitely block represented in every block
subspace of $T[(\FA_k,\theta_k)_{k=1}^\infty]$.

\dem Given $n \in \N$ and $0 < \varepsilon < 1$, take $k \in \N$
such that $n_k > n$ and $\ds n_k
-\left(\lambda_{n_k^{\ell_k}}\right)^{1 \over \ell_k} <
\varepsilon$.  Let $(u_i)_{i=1}^\infty$ be a normalized block
sequence. Then $$ n_k^{\ell_k} \geq \norma
{\sum_{i=1}^{n_k^{\ell_k}} u_i} \geq \norma
{\sum_{i=1}^{n_k^{\ell_k}} e_i} = \lambda_{n_k^{\ell_k}} \geq
(n_k-\varepsilon)^{\ell_k} $$ and, by Proposition 8,
$\ell_1^{n_k}$ is finitely block represented in blocks of
$(u_i)_{i=1}^\infty$. \bull

\proclaim Remark 5. By similar arguments it is easy to prove that
the following condition is also sufficient:
\item {1.} There exits $m \geq 2$ such that
$\ds \lim_{\ell \to \infty} (\lambda_{m^\ell})^{1 \over \ell} =
m$.

\noindent We can also give sufficient conditions for the sequence
$(\theta_k)_{k=1}^\infty$:
\item {2.} There exists $(n_k)_{k=1}^\infty \subseteq \N$ unbounded
 and $(\ell_k)_{k=1}^\infty$ such that
$\ds \lim_{k\to \infty} n_k \left[1 -
\left(\theta_{n_k^{\ell_k}}\right)^{1 \over \ell_k} \right] = 0 $
or
\item {3.} There exists $m \geq 2$ such that
$\ds
\lim_{\ell \to \infty} (\theta_{m^\ell})^{1 \over \ell} =1
$
or, equivalently, $\ds \lim_{\ell \to \infty} (\theta_{m^\ell})^{1
\over \ell} =1$ for all $m \geq 2$.

\proclaim Lemma 5. Let $(X, \norma\cdot)$ and $(X',
{\norma\cdot}')$
 be Banach spaces not
totally incomparable with Schauder bases $(e_i)_{i=1}^\infty$ and
$(e_i')_{i=1}^\infty$. If $(e_i)_{i=1}^\infty$ is shrinking, there
exist  block sequences $(u_i)_{i=1}^\infty$ and
$(u_i')_{i=1}^\infty$ of $(e_i)_{i=1}^\infty$ and
$(e_i')_{i=1}^\infty$ respectively such that the application $T:
\overline {Span} \{ u_i\ \vert\ i \in \N\} \lra \overline {Span}
\{ u_i'\ \vert\ i \in \N\}$, given by $T(u_i)=u_i'$ for all $i \in
\N$ is an isomorphism.

\dem  There exist subspaces $Y \subseteq X$ and $Y' \subseteq
X'$  and an isomorphism $S: Y \lra Y'$. We will see that for all
$\varepsilon>0$ we can find block sequences $(u_i)_{i=1}^\infty$
and $(u_i')_{i=1}^\infty$ such that $\ds (1-\varepsilon){\norma
S} \Vert{S^{-1}}\Vert \le {\norma T} \Vert{ T^{-1}}\Vert\le
(1+\varepsilon){\norma S} \Vert{S^{-1}}\Vert$.

Let $\varepsilon>0$. There exists a normalized block sequence
$(x_i)_{i=1}^\infty$ of
 $(e_i)_{i=1}^\infty$ and $\overline {Span} \{ y_i\ \vert\ i\in \N \}
 \subseteq Y$ such that the linear isomorphism defined by
 $U(x_i)=y_i$ verifies ${\norma U } \Vert{U^{-1}}\Vert\le
 1+\varepsilon$. Let $y_i':=S(y_i)$ for all $i \in \N$.

Since $\inf_{i\in\N}\norma{y_i'}>0$ and $(e_i)_{i=1}^\infty$ is a
 shrinking basis,  $y_i'$ tends to $0$ weakly.
So, by the Bessaga-Pelcynski principle, there is a subsequence
  $(y_{i_k}')_{k=1}^\infty$  and a block sequence
$(u_k')_{k=1}^\infty$ of $(e_i')_{i=1}^\infty$ such that the
isomorphism defined by $V(y_{i_k}')=u_k'$ verifies ${\norma V}
\Vert{V^{-1}}\Vert\le
 1+\varepsilon$. Take $u_k=x_{i_k}$ and $T=V\circ S\circ U$. \bull

\proclaim Remark 6. Let $X=T[(\FA_k,\theta_k)_{k=1}^\infty],
\theta_k\in(0,1)$. Since its canonical basis $(e_i)_{i=1}^\infty$
is unconditional, hence being shrinking is equivalent to $\ell_1$
not being isomorphic to any subspace of $X$ and this is the case
by Proposition 7.

\proclaim Theorem 4. Let $X=T[(\FA_k,\theta_k)_{k=1}^\infty]$ and
$X'=T[(\FA_k,\theta_k')_{k=1}^\infty]$ with $\theta_k,
\theta_{k'}\in(0,1)$ be
 such that
$\ell_1$ is finitely block represented in every block subspace of
$X$ and $X'$. If $X$ and $X'$ are not totally incomparable, then
there exists $C \ge 0$ such that for all $n \in \N$, $$ {1 \over
C} \leq {\lambda_\ell \over \lambda_\ell'} \leq C. \leqno {(*)} $$

\dem Denote by $\norma\cdot$ and ${\norma\cdot}'$ the norms of $X$
and $X'$, respectively. By Lemma 5, there exist block sequences
$(u_i)_{i=1}^\infty\subseteq X$ and $(u_i')_{i=1}^\infty \subseteq
X'$  of their respective bases denoted by $(e_i)_{i=1}^\infty$ and
$(e_i')_{i=1}^\infty$, such that $T: \overline {Span} \{ u_i\
\vert\ i \in \N\} \lra \overline {Span} \{ u_i'\ \vert\ i \in
\N\}$, given by $T(u_i)=u_i'$ for all $i \in \N$ is an
isomorphism. Therefore, for all $(a_i)_{i=1}^\infty \subseteq \R$
and $n \in \N$,
$$ {1 \over {\norma T}} \norma {\sum_{i=1}^n a_i
u_i'}'\leq \norma {\sum_{i=1}^n a_i u_i} \leq {\norma {T^{-1}}}\
\norma {\sum_{i=1}^n a_i u_i'}'. $$

By Lemma 4, given $\varepsilon
>0$ and $\ell \in \N$, there exists a
normalized block sequence $y_1,\dots, y_{\ell}$  of
$(u_i)_{i=1}^\infty$, such that $$ \lambda_\ell-\varepsilon \leq
\norma {\sum_{i=1}^\ell y_{i}} \leq \lambda_\ell +\varepsilon. $$

Let $y_{i}':=T(y_{i})$ for all $i=1,\dots,\ell$. Then we have $$
\eqalign { \lambda_\ell +\varepsilon &\geq \norma {\sum_{i=1}^\ell
y_{i}} \geq {1 \over {\norma T}} \norma {\sum_{i=1}^\ell y_{i}'}'
=\cr &={1 \over {\norma T}} \norma {\sum_{i=1}^\ell
\norma{y_{i}'}' {y_{i}' \over \norma{y_{i}'}'}}' \geq {1 \over
{\norma T}} \min_{1 \leq i \leq \ell} \norma {y_{i}'}' \norma
{\sum_{i=1}^\ell {y_{i}' \over \norma{y_{i}'}'}}' \geq \cr &\geq
{1 \over {\norma T} {\norma {T^{-1}}}} \norma {\sum_{i=1}^\ell
e_i'}'= {1 \over {\norma T} {\norma {T^{-1}}}} \lambda_\ell' \cr}
$$
(note that in the last inequality we use Lemma 3). Since the
inequality is true for all $\varepsilon >0$, we have proved that
$\ds \lambda_\ell \geq {1 \over {\norma T}{\norma {T^{-1}}}}
 \lambda_\ell'.$

Now we reverse the roles of $X$ and $X'$ to obtain $\ds {1 \over
{\norma T} {\norma {T^{-1}}}} \lambda_\ell' \leq \lambda_\ell \leq
{\norma T} \Vert T^{-1}\Vert \lambda_\ell'. $ \bull

\proclaim Remark 7. If $X$ and $X'$ contain isometric subspaces
$Y$ and $Y'$ , then $\lambda_\ell = \lambda_\ell'$ for all $\ell
\in \N$. Actually, the same equality holds if for every
$\varepsilon > 0$, $X$ and $X'$ contain
$(1+\varepsilon)-$isomorphic subspaces.

\proclaim Remark 8. There are special cases when the calculus of
$\lambda_\ell$ is easy. For instance when $(\theta_k),
(\theta'_k) $ belong to the so called class $\cal F$ defined in
[11] we have $ \lambda_\ell= \ell \cdot \theta_\ell $ and the
condition $(*)$ of Theorem 4 yields $\ds {1 \over C} \leq
{\theta_\ell \over \theta'_\ell} \leq C $ for all $\ell$ or
$\theta_\ell=\theta_\ell'$ if we can find isometric subspaces or
$(1+\varepsilon)-$isomorphic subspaces for all $\varepsilon>0.$

\proclaim Example 3. Let $f_r(x)=\log_2^{r}(1+x)$ with $0 < r < 3
\log 2-1$. Then $(f_r^{-1}(k))\in\cal F$ and if $0 < r < s <
3\log 2 - 1$, the spaces $\ds T\left[\left(\FA_k, {1\over f_r(k)}
\right)_{k=1}^\infty \right]$ and $\ds T\left[\left(\FA_k, {1
\over f_s(k)} \right)_{k=1}^\infty \right]$ are, by Theorem 4,
totally incomparable. Moreover, it is easy to check that these
spaces are also totally incomparable to $\ell_p$, $1 \leq p <
\infty$ or $c_0$.

\parrafo{References}
\item{[1]} S. A. Argyros and I. Deliyanni,
{\it Banach spaces of the type of Tsirelson},
Preprint,1992.
\item{[2]} S. A. Argyros and I. Deliyanni,
{\it Examples of asymptotic $\ell^1$ Banach spaces}, Trans. Amer.
Math. Soc, {\bf 349} (1997), 973-995.
\item{[3]} Bellenot, {\it Tsirelson superspaces and $\ell_p$},
Journal of Functional Analysis, {\bf 69} (1986) 207-228.
\item{[4]} J. Bernu\'es and I. Deliyanni,
{\it Families of finite subsets of $\N$ of low complexity and Tsirelson
type spaces},
To appear in Math. Nach.
\item{[5]} J. Bernu\'es and Th. Schlumprecht,
{\it El problema de la distorsi\'on y el problema de la base
incondicional},
Colloquium del departamento de an\'alisis. Universidad Complutense.
Secci\'on 1, {\bf 33} (1995).
\item{[6]} P.G. Casazza and T. Shura,
{\it Tsirelson's Space}, LNM 1363, Springer-Verlag, Berl\'{\i}n, 1989.
\item{[7]} T. Figiel and W.B. Johnson,
{\it A uniformly convex Banach space which contains no $\ell_p$},
Compositio Math. {\bf 29} (1974), 179-190.
\item{[8]} J. Lindenstrauss and L. Tzafriri,
{\it Classical Banach Spaces I, II}, Springer-Verlag, New York,
1977.
\item{[9]} A. Manoussakis, {\it On the structure of a certain class of mixed
Tsirelson spaces}, Preprint, 1999.
\item{[10]} E. Odell and T. Schlumprecht,
{\it A Banach space block finitely universal for monotone basis}.
To appear in  Transactions of the Amer. Math. Soc.
\item{[11]} Th. Schlumprecht,
{\it An arbitrarily distortable Banach space},
Israel Journal of Math. {\bf 76} (1991), 81-95.
\item{[12]} B.S. Tsirelson, {\it Not every Banach space contains an embedding
of $\ell_p$ or $c_0$}, Funct. Analysis and Appl. {\bf 8} (1974),
138-141.
\item{[13]} L. Tzafriri, {\it On the type and cotype of Banach spaces},
Israel Journal of Math. {\bf 32} (1979), 32-38.
\bigskip\bigskip\bigskip

Julio Bernu\'{e}s and Javier Pascual\bigskip

Departamento de Matem\'{a}ticas

Universidad de Zaragoza.

50009-Zaragoza (Espa\~{n}a)\bigskip

e-mail: bernues@posta.unizar.es
\end